\documentclass[12pt]{article}
\usepackage{amsfonts,amsmath,amssymb,amscd,amsthm}
\usepackage{graphicx}
\long\def\comment#1\endcomment{}

\begin{document}

\newpage
\centerline
{\uppercase
{\large {\bf
 A smaller counterexample to the Lando conjecture}} \footnote{This paper is prepared under the supervision of Arkadiy Skopenkov and is submitted to the Moscow Mathematical Conference for High-School Students. Readers are invited to send their remarks and reports on this paper to mmks@mccme.ru}}

\bigskip
\centerline{\bf  V.Belousov }

\bigskip

{\bf Abstract. } The following conjecture was proposed in 2010 by S. Lando.

Let $M$ and $N$ be two unions of the same number of disjoint circles in a sphere. Then there exist two spheres in 3-space whose intersection is transversal and is a union of disjoint circles that is situated as $M$ in one sphere and as $N$ in the other . Define union $M$ of disjoint circles to be {\it situated} in one sphere as union $M_1$ of disjoint circles in the other sphere if there is a homeomorphism between these two spheres which maps $M$ to $M_1$.
%Avvakumov found the first counterexample to Lando's conjecture. Then he found the minimal counterexample. 

%Avvakumov found special $M_1$ and $N_1$, each of which is a union of 9 circles in a sphere. And using Avvakumov theorem and brute force computer method he proved that there exists no pair of spheres in 3-space whose intersection is transversal and is a union of disjoint circles that is situated as $M_1$ in one sphere and as $N_1$ in the other. That was the first counterexample to Lando's Conjecture. Then he found special $M_2$ and $N_2$, each of which is a union of 7 circles in a sphere. The aim of this work is prove of that there exists no pair of spheres in 3-space whose intersection is transversal and is a union of disjoint circles that is situated as $M_2$ in one sphere and as $N_2$ in the other. Pair $M_2$, $N_2$ is the minimal counterexaple to Lando's Conjecture.
In this paper we prove that there exists pair of sets of 7 circles in sphere, that is a counterexample to the Lando conjecture. This is proved using the Avvakumov Theorem. 
%Also if pair $(M_1,N_1)$ is a counterexample to the Lando conjecture then $M_1$ and $N_1$ contain at least 7 circles.
We conjecture that there exists no pair $(M,N)$ that is counterexample and $M$ contains 6 or less circles.

%$W$ is a certain set of  $7$ disjoint circles in sphere. 
%$R$ is another certain set of 7 disjoint  circles in sphere. The aim of this paper is proof of  that there doesn't exist  pair of topological sphere in 3-space, such that their intersection is $W$ in the first sphere and $R$ in another. The previous proof was brute force %computer method. Both proofs  use Avvakumov theorem, that shows the relationship between a topological problem and a graph problem.  

 \bigskip
{\bf Definitions.}

Let $p$ and $q$ be two sets of edges of a tree $Y$. 

The set $p$ is {\it  on the same side} of $q$ (in this tree $Y$) if $p \bigcap q = \emptyset$ and for each two vertices of edges of $p$ there is a path in the tree connecting these two vertices,
%!!
  and containing an even number of edges of $q$. Sets $p$ and $q$ are {\it unlinked} (in this tree) if $p$ is on the same side of $q$ and $q$ is on the same side of $p$ .

For vertice $P$ of graph we denote as $\delta P$ all edges whose end is $P$.

Let $K$ and $K'$ be two trees with the same number of edges.
Let $h$ be a bijection (i.e. one-to-one correspondence) between their edges.

 Then $h$ is called {\it realizable} if $h(\delta A)$ and $h(\delta B)$ are unlinked for each two vertices $A$ and $B$ in $K$ such that the path joining $A$ and $B$ contains even number of edges.

Graphs $K$ and $K'$ are {\it friendly} if such a bijection exists.
%$K$ is a friend of  $K'$  if such bijection exists.
%If $K$ is a friend of  $K'$  and $K'$ is a friend of $K$ then $K$ and $K'$ are friendly.

 %Note: Avvakomov Theorem implies that it is sufficient $K$ to be a friend of  $K'$ .

\smallskip

 	\begin{figure}[h]
	\centering\includegraphics[width=15cm]{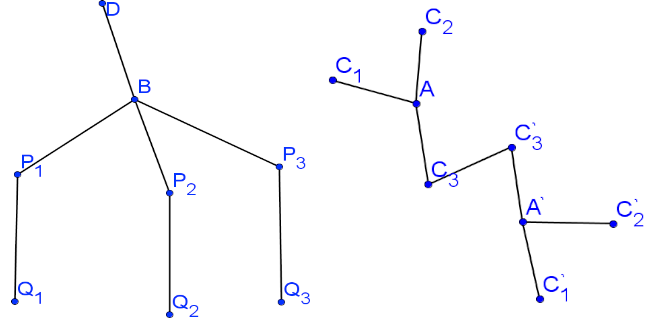}
 	\caption{Graphs $H$ and $G$.}
%	\label{f-trans}
	\end{figure}
\smallskip

%Graph $G$ has vertices $A$,$C_3$,$C_1$,$C_2$,$A'$,$C'_3$,$C'_1$,$C'_2$ and edges $C_3C'_3$, $AC_i$, $A'C'_i$, $i = 1,2,3$.
%Graph $H$ has vertices $B$,$D$,$P_1$,$P_2$,$P_3$,$Q_1$,$Q_2$,$Q_3$ and edges $BD$,$BP_i$,$P_iQ_i$, $i = 1,2,3$.
Let graph $G$ be a graph that has vertices $A$,$C_i$,$A'$,$C'_i$ and edges $C_3C'_3$, $AC_i$, $A'C'_i$, $i = 1,2,3$.
Let graph $H$ be a graph that has vertices $B$,$D$,$P_i$,$Q_i$ and edges $BD$,$BP_i$,$P_iQ_i$, $i = 1,2,3$.

\bigskip
{\bf Theorem 1. }{\it Graphs $G$ and $H$ are unfriendly. }
 \smallskip
 
  Lets state a result that shows why this theorem is interesting. 
 Suppose that $M$ is a union of disjoint circles in sphere $S^2$. Define (`dual to $M$') graph $G = G(S^2;M)$ as follows. The vertices are
the connected components of $S^2$ \textbackslash
 $M$. Two vertices are connected by an edge if the corresponding connected components are neighbors.

\smallskip
 {\bf Avvakumov Theorem. }[A]
{\it
%Suppose that $M$ is a union of disjoint circles in sphere $S^2$.Suppose that $N$ is a union of disjoint circles in sphere $S^2$.
 Let $M$ and $N$ be two unions of the same number of disjoint circles in a sphere $S^2$.
 Then there exist two spheres in 3-space whose intersection is
transversal and is a union of disjoint circles that is situated as $M$ in one sphere 
%(union $X$of disjoint circles is situated in one sphere as union $Y$ of disjoint circles in the other sphere if there is a homeomorphism between these two spheres which maps $X$ to $Y$)
 and as $N$ in the other if and only if the graph dual to $M$  and $N$ are friendly.
 }

% You can see prove of this theorem in [A].
 %Note: Avvakumov theorem implies that if K and K' are friendly, then K' and K are friendly.
 This theorem implies that  friendliness is symmetric. This will be used in the proof.
% $G$ and $H$ are dual to $M_2$ and $N_2$.

\bigskip
%{\it Proof of theorem 1.}
%Then there is bijection between their edges with the above properties.
%By $H(n)$ we denote the image in graph $H$ of edge $n$ of graph $G$,
%and by $G(n)$ we denote the image in graph $G$ of edge $n$ of graph $H$ (bijection inversed). Let $X$ be a subset of the set of edges of graph $H$.
%$G(X)$ is graph, formed by images of edges from $X$.{

Suppose $\phi$ is a realizable bijection between edges of $G$ and $H$. For edges $e_1,e_2,e_3,\ldots ,e_n$ of graph $G$ by $h (e_1,e_2,e_3,\ldots,e_n)$ we denote subgraph formed by $\phi (e_1),\phi (e_2),\phi (e_3),\ldots,\phi(e_n)$ in graph $H$.
And for edges $e_1,e_2,e_3,\ldots,e_n$ of graph $H$ by $g(e_1,e_2,e_3,\ldots,e_n)$ we denote subgraph formed by  $\phi(e_1),\phi(e_2),\phi(e_3),\ldots ,\phi(e_n)$ in graph $G$.
%(bijection inversed).
% g и h назвать отображения. 
% представление работ на ммкш.выыыв

\bigskip
{\bf Proposition 1.}
{\it Both graphs $H_1:=h(AC_1,AC_2,AC_3)$ and $H_2:=h(A'C'_1,A'C'_2,A'C'_3)$ are connected.}

\smallskip

{\it Proof.} Lets prove the connectedness for $H_1$, and for $H_2$ the proof is analogous.

If $H_1$ is not connected then one of edges from $H\setminus H_1 = h(AC'_1,AC'_2,AC'_3,C_3C_3')$ belongs to path connecting two edges from $H_1$.
Vertices $A$ and $C_1'$ are linked by a path of even length.
So $h(A'C'_1)$ doesn't belong to any path that joins a pair of edges of graph $H_1$. Analogically $h(A'C'_2)$ doesn't belong to any path that joins a pair of edges of graph $H_1$.
Vertices $A$ and $C'_3$ are linked by a path of even length too. Hence,

$\bullet$ {\it Case 1.} Neither $h(C_3C'_3)$ nor $h(A'C'_3)$ don't belong to any path that joins a pair of edges of graph $H_1$; 

$\bullet$ {\it Case 2.} $h(C_3C'_3)$ and $h(A'C'_3)$ belong to path that joins a pair of edges $J_1$,$J_2$ of graph $H_1$.

In the first case the graph $H_1$ is connected.

In the second case $J_1$,$h(C_3C'_3)$,$h(A'C'_3)$,$J_2$ form a path of length 4. Without loss of generality this path is $Q_1P_1BP_2Q_2$.
Hence path, that links edges $J_1$ and $H_1-J_1-J_2$, intersect only one of edges $h(C_3C'_3)$,$h(A'C'_3)$.
Which is impossible. 

QED

\bigskip

{\bf Proposition 2.}
{\it Vertex $B$ is an endpoint of edge $h(C_3C'_3)$. }
\smallskip

{\it Proof.}
 Vertices $B$ and $Q_i$ are linked by a path of even length. $g(BD,BP_1,BP_2,BP_3) = g(\delta B)$ is unlinked with any edge $g(P_iQ_i)$ in $G$. This implies $g(\delta B)$ is connected.
%Hence $G(BD,BP_1,BP_2,BP_3) = G(\delta B)$ is connected, because it is unlinked with any edge of $G(P_iQ_i)$ in $G$ .
There are only 2 connected subgraphs with 4 edges in $G$ up to automorphism of $G$:

The first subgraph, say $X$, has vertices $C_1,C_2,C_3,A,C'_3$ and its edges are precisely the edges of $G$ with both ends in $X$.

The second subgraph, say $Y$, has vertices are $C_1,A,C_3,C'_3,A'$ and its edges are precisely the edges of $G$ with both ends in $Y$.

Since $C_3C_3'$ is fixed under $\mathrm{Aut} (G)$  and  $C_3C_3'$ is contained both in $X$ and $Y$, one of edges $BD,BP_1,BP_2,BP_3$ is $h(C_3C'_3)$.

QED

\bigskip

{\it Proof of theorem 1.}
Suppose graphs $G$ and $H$ are friendly. Then there exists a realizable bijection $\phi$ between edges of $G$ and $H$.

According to proposition 1 graph $H-h(C_3C'_3)$ is a union of two connected
graphs with 3 edges.
 %Hence one of them containes at least two of the edges $P_iQ_i$ as his edge.
  According to proposition 2 one of these graphs containes at least two of the edges $P_iQ_i$ as his edge. 
Without loss of generality let $P_1Q_1$ and $P_2Q_2$ be in $h(AC_1,AC_2,AC_3) = H_1$.
But then the length of the path linking $Q_1$, $Q_2$ is 4. Since there are only 3 edges in $H_1$, this is impossible.

 QED
 
\bigskip
\bigskip

{\bf  References}
\smallskip
%добавить ссылка на ссылки

[A]S. Avvakumov, A counterexample to the Lando conjecture on intersection of spheres in 3-space,
http://arxiv.org/pdf/1210.7361v2.pdf

For related resultes see:

[B]A. Rukhovich, On intersection of two embedded spheres in 3-space, http://arxiv.org/abs/ 1012.0925

[C]S. Avvakumov, A. Berdnikov, A. Rukhovich and A. Skopenkov, How do curved spheres intersect in 3-space, or two-dimensial meandra, http://www.turgor.ru/lktg/2012/3/3-1en\_si.pdf

[D]T. Hirasa, Dissecting the torus by immersions, Geometriae Dedicata, 145:1 (2010), 33-41

[E]T. Nowik, Dissecting the 2-sphere by immersions, Geometriae Dedicata 127, (2007), 37-41,
 http: //arxiv.org/abs/math/0612796

\end{document}